\documentclass[8pt]{article}

\usepackage{flushend, cuted}
\usepackage{amsmath}
\usepackage{amsfonts}
\usepackage{bbm}

\usepackage{amssymb}
\usepackage{float}
\usepackage{color,amsmath,amssymb, amsfonts,amstext,amsthm}

\usepackage{epsfig, graphicx, graphics}
\usepackage{longtable}
\usepackage{cite}
\usepackage{subfigure}

\renewcommand{\phi}{\varphi}

\title{Geometric Methods for Stochastic Dynamical Systems \footnote{
$\quad$$^{1}$ Department of Applied Mathematics, Illinois Institute of Technology, Chicago 60616, USA.
$\quad$$^{2}$ School of Mathematics and Statistics, Zhengzhou University, 100 Kexue Road, Zhengzhou 450001, China.}}

\author{Jinqiao Duan$^{1}$, Hui Wang$^{2}$}

\date{\today}

\begin{document}

\maketitle

\section{Introduction}\label{sec1.1}

Noisy fluctuations are ubiquitous in complex systems. They play a crucial or delicate role in the dynamical evolution of  gene regulation, signal transduction, biochemical reactions, among other systems. Therefore, it is essential to consider the effects of noise on   dynamical systems. It has been a challenging topic to  have  better understanding of the impact of the noise on the dynamical behaviors of complex  systems. See \cite{Arnold2003, Duan2015} for more information.

\subsection{Stochastic dynamical systems}

 A dynamical system may be thought as evolution of mechanical `particles',  and is then described by a system of ordinary differential equations
\begin{equation}
\frac{dx}{dt}=f(x),
\end{equation}
where $f$ is often called a vector field.
In natural science and applied science, the dynamical systems are used to describe the
evolution of complex phenomena \cite{Meiss2017, Wiggins2003, Guckenheimer1983}.  However,  dynamical systems are often  influenced by   random factors in the environment, such as in the systems of the propagation of waves through random media, stochastic particle acceleration, signal detection, and optimal control with fluctuating constraints. Indeed, a small random disturbance may have an unexpected  effect on the whole dynamical system.  A stochastic dynamical system is a  dynamical system with noisy components or under random influences  \cite{Arnold2003, Duan2015} and thus may be modeled by a stochastic differential equation
\begin{equation}
\frac{dx}{dt}=f(x) + \mbox{noise}.
\end{equation}

Usually, we use two kinds of   noise. One is  Gaussian noise, and  the other is non-Gaussian noise.  Gaussian noise is modeled by Brownian motion $B_t$, while  non-Gaussian noise is expressed via L\'evy process (especially $\alpha$-stable L\'evy motion  $L_t^{\alpha, \beta}$ )  \cite{Applebaum2009, Sato1999}.  Dynamical systems driven by Gaussian noise have been widely studied, but in some complex systems, the random influences or stochastic processes are non-Gaussian. For example,  during the regulation of gene expression,   transcriptions of DNA from genes and translations into proteins  take  place in a bursty, intermittent, unpredictable manner \cite{Holloway2017, Kumar2015, Dar2012, Ozbudak2002, Blake2003, Sanchez2013}. It is more suitable to  model these processes  by dynamical systems with non-Gaussian L\'evy noise.

In order to better describe the impact of noise on the dynamical systems, we consider some  geometric methods. These include most probable phase portraits, mean phase portraits,  invariant manifolds, and slow manifolds. These geometric tools may help us characterize the impact of noise on the dynamical systems vividly from the geometric perspective.

\section{Stochastic dynamial systems}

Consider a stochastic dynamical system in $n$-dimensional Euclidean space $\mathbb{R}^n$, either with (Gaussian) Brownian noise
\begin{equation}\label{s1}
dX_t=f(X_t)dt + \sigma(X_t)B_t, \;\;  X_0=x_0,
\end{equation}
or (non-Gaussian) L\'evy noise
\begin{equation}\label{s11}
dX_t=f(X_t)dt + \sigma(X_t) L_t^{\alpha, \beta},   \;\; X_0=x_0,
\end{equation}
where $\sigma$ is noise intensity.
Let us briefly recall the definition of Brownian motion and L\'evy motion.

\subsection{Brownian motion}

\textbf{Definition}(\cite{Karatzas1991}) A Brownian motion $B_t$ is a stochastic process defined on a probability space $\Omega$
equipped with probability $\mathbb{P}$,  with the following properties:\\
(i)  $B_{0} = 0$,  almost surely;\\
(ii) $B_t$ has independent increments;\\
(iii) $B_t$ has  stationary increments with normal distribution: $B_t-B_s  \sim  \mathcal{N}(t-s, 0)$,  for $t>s$;\\
(iv) $B_t$ has  continuous sample paths, almost surely.

\subsection{L\'evy motion}
\textbf{Definition}(\cite{Duan2015})   On a sample space $\Omega$ equipped with probability $\mathbb{P}$, a scalar asymmetric   stable L\'{e}vy motion $L_t^{\alpha, \beta}$,  with the non-Gaussianity index $\alpha \in (0,2)$ and the skewness index $\beta \in [-1,1]$,   is a stochastic process with the following properties:\\
(i)  $L_{0}^{\alpha, \beta} = 0$,  almost surely; \\
(ii) $L_t^{\alpha, \beta}$ has independent increments;\\
(iii) $L_t^{\alpha, \beta}$ has stationary increments with stable distribution: $L_t^{\alpha, \beta}-L_s^{\alpha, \beta} \sim $ $S_\alpha((t-s)^\frac{1}{\alpha}, \beta, 0)$, for $t\textgreater s$;\\
(iv) $L_t^{\alpha, \beta}$ has stochastically continuous sample paths, i.e., for every $s$, $L_t^{\alpha, \beta} \rightarrow $ $L_s^{\alpha, \beta}$ in probability (i.e., for all $\delta>0$, $\mathbb{P}$($\mid$$L_t^{\alpha, \beta}-L_s^{\alpha, \beta}$$\mid$ $\textgreater$ $\delta$)$\rightarrow$ 0 ),  as $t\rightarrow s$.

The jump measure,  which describes jump intensity and size for sample paths,  for the asymmetric L\'{e}vy motion  $L_t^{\alpha, \beta}$  is \cite{Applebaum2009,Duan2015},
\begin{equation} \label{eq:3}
 \nu_{\alpha, \beta}(dy)=\frac{C_1 I_ {\{0<y<\infty\}}(y)+C_2 I_{\{-\infty<y<0\}}(y)}{\mid y\mid ^{1+\alpha}}dy ,
\end{equation}
with $ C_1 =\frac{H_\alpha (1+\beta)}{2}$ and $C_2= \frac{H_\alpha (1-\beta)}{2}.$
When $\alpha=1$,  $H_\alpha=\frac{2}{\pi}$; when $\alpha\neq 1$, $H_\alpha=\frac{\alpha(1-\alpha)}{\Gamma(2-\alpha)\cos(\frac{\pi \alpha}{2})}.$\\

Especially for  $\beta=0$,  this is the     symmetric stable L\'{e}vy motion,  which is usually denoted by  $L_t^\alpha  \triangleq  L_t^{\alpha, 0}$. More especially, for $\alpha=2, \beta=0$, this is the well-known Brownian motion $B_t$.

\section{Geometric methods for stochastic dynamical systems}\label{sec2}
As in the geometrical approaches for deterministic dynamical systems\cite{Meiss2017, Arnold1988, Arnold1989},  we consider   geometric  concepts and methods for stochastic dynamical systems qualitatively. These include
phase portraits and invariant manifolds.

We now introduce some geometric methods for stochastic dynamical systems in this section.

\subsection{Most probable phase portraits}

The Fokker-Planck  equation for  the stochastic differential equation \eqref{s1} or \eqref{s11} describes the time evolution of  the probability density $p(x,t)   \triangleq p(x, t|x_0, 0)$ for the solution process $X_t $ with initial condition $X_0=x_0$. It is a linear deterministic partial differential equation \cite{Duan2015}
\begin{equation} \label{fpe}
 p_t  = A^* p,  \;\;   p(x,0)=\delta(x-x_0),
\end{equation}
where  $A^*$ is the adjoint operator of the generator  $A$ for this stochastic differential equation, and $\delta$ is the Dirac delta function.

As the solution of the Fokker-Planck equation,  the probability density function $p(x,t)$ is a surface in the $(x,t,p)$-space. At a given time instant $t$, the maximizer $x_m(t)$ for $p(x,t)$ indicates the most probable (i.e., maximal likely) location of this orbit at time $t$. The orbit traced out by $x_m(t)$ is called a most probable orbit  starting at $x_0$. Thus,  the deterministic orbit $x_m(t)$ -- also denoted by $x_m(x_0, t)$  --- follows the top ridge   of the surface in the $(x,t,p)$-space as time goes on.

\textbf{Definition}{(Most probable equilibrium point)}
A most probable  equilibrium point (state) is a  point (state) which either attracts or repels  all nearby points (states). When it attracts  all nearby points (states), it is called a   most probable \emph{stable} equilibrium point (state), while if it repels all nearby points (states), it is called a   most probable \emph{unstable} equilibrium point (state).


The most probable phase portrait \cite{Duan2015,Cheng2016}  for a stochastic dynamical system  is the state space with representative most probable orbits including  equilibrium states. It is a   deterministic  geometric object, which describes which sample orbit is most probable (maximal likely).

\subsection{Mean phase portraits}
Given the probability density function $p(x,t)   \triangleq p(x, t|x_0, 0) $ as the solution of the Fokker-Planck equation \eqref{fpe}. The mean orbit starting at an initial point $x_0$ in state space is defined as
\begin{equation}
 \bar{x}(x_0, t)=\int_\mathbb{R} \xi  \; p(\xi,t|x_0,0) d\xi.
\end{equation}

We also give the definition of a mean equilibrium point (state).

\textbf{Definition}{(Mean equilibrium point)}
A mean equilibrium point (state) is a  point (state) which either attracts or repels  all nearby points (states). When it attracts  all nearby points (states), it is called a   mean \emph{stable} equilibrium point (state), while if it repels all nearby points (states), it is called a   mean \emph{unstable} equilibrium point (state).

The mean probable phase portrait is the state space with representative mean orbits including mean equilibrium states.   It is also a deterministic  geometric object, which describes the   orbits in the mean sense,  starting  in state space.

\subsection{Invariant manifolds}
Invariant geometry structures play an important role in our understanding of dynamical systems. Invariant manifolds to nonlinear systems just as eigenspace to linear systems.  With a better understanding of stochastic invariant manifolds, we may have a new insight to understand stochastic dynamical systems.

 We recall some definitions.

\textbf{Definition}{(Random set)}
A random set for a random dynamical $\varphi$ in $\mathbb{R}^n$ is a collection $M=M(\omega)$, $\omega \in \Omega $, satisfies the following property:\\
(i) $M(\omega)$ is a nonempty closed set, $M(\omega)\subset \mathbb{R}^n $, $\forall $ $\omega \in \Omega$;  \\
(ii)$V_x: \Omega \rightarrow \mathbb{R}^1  $, and $V_x(\omega)\triangleq inf_{y\in M(\omega)}d(x,y)$ is a scale random variable $\forall$ $x\in \mathbb{R}^n$.

\textbf{Definition}{(Random invariant set)}
An invariant set for a random dynamical $\varphi$ is a random set $M$ satisfies:
$\varphi(t,\omega,M(\omega))=M(\theta_t \omega)$, $\forall$ $t\in \mathbb{R}$ and $\omega \in \Omega$.

\textbf{Definition}{(Random invariant manifold)}
If a random invariant set $M$ for a random dynamical $\varphi$ can be represented by a graph of a Lipschitz mapping:
$\gamma^*(\omega, \cdot): H^+\rightarrow H^-$, with direct sum decomposition $ H^ \oplus H^-=\mathbb{R}^n$, satisfies that $M(\omega)={(x^+,\gamma^*(\omega, x^+ )), x^+\in H^+}$, then $M$ is called a Lipschitz invariant manifold.

For more information about random invariant manifolds, see \cite{Arnold2003, Duan2003, Fu2008}.

\subsection{Slow manifolds}

When a system has a slow component and a fast component, we call it slow-fast system. For example, temperature and molecular motion, climate and weather,  are slow and fast components in a body and in the climate system, respectively. In this kind of systems, there are   two time scales, and usually,  we are more interested in the slow dynamics. But we still like to take the fast dynamics into account, even when the fast evolution is less important to our modeling purpose.  We thus represent fast components in terms of slow components as slow manifolds.    For more information about slow manifolds, see \cite{Fu2008, Fu2013, Ren2015a, Ren2015b, Constable2013}.

Consider a stochastic slow-fast system:
\begin{eqnarray}\label{eqs}
\begin{cases}
   \dot{x} & = Ax + f(x,y), x(0) = x_0\in \mathbb{R}^n,\\
   \dot{y} & = \frac{1}{\varepsilon} By+\frac{1}{\varepsilon}g(x,y)+\frac{\sigma}{\sqrt {\varepsilon} \dot{W_t}},
   \end{cases}
\end{eqnarray}
where $A$ and $B$ are matrices, $\varepsilon$ is a small positive parameter measuring slow and fast scale separation, $f$ and $g$ are nonlinear Lipschitz continuous functions with Lipschitz constant $L_f$ and $L_g$ respectively, $\sigma$ is a noise intensity constant, and ${W_t : t \in \mathbb{R}}$ is a two-sided $\mathbb{R}^m-valued$ Brownian motion on a probability space $(\Omega, \mathcal{F},\mathbb{P})$.

Under the exponential dichotomy and gap condition, there exists a random slow manifold $\tilde{\mathcal{M}}^\varepsilon (\omega)={(\xi, \tilde{h}^\varepsilon (\xi,\omega)):\xi\in \mathbb{R}^n}$ with $\tilde{h}^\varepsilon$  can be expressed by
Liapunov-Perron equation:
\begin{equation}
\tilde{h}^\varepsilon (\xi,\omega)=\frac{1}{\varepsilon}\int_{-\infty}^0 e^{-\frac{B}{\varepsilon}s}g(X(s,\omega,\xi),Y(s,\omega,\xi)+\sigma\eta^\varepsilon(\theta_s \omega))ds
\end{equation}
Then the slow system is
\begin{equation}
\dot{x}=Ax+f(x,\tilde{h}^\varepsilon(x,\theta_t \omega)+\sigma\eta(\psi_\varepsilon \omega))
\end{equation}
Moreover, $\hat{h}^\varepsilon(\xi,\omega)=\tilde{h}_d^\varepsilon(\xi,\omega)+\varepsilon\tilde{h}_1^\varepsilon(\xi,\omega)$
  is an approximation or a first order truncation of  $\tilde{h}^\varepsilon(\xi,\omega)$. Thus,  we have an approximate slow system
  \begin{equation}
  \dot{x}=Ax+f(x,\hat{h}^\varepsilon(x,\theta_t \omega)+\sigma\eta(\psi_\varepsilon \omega))
  \end{equation}
 Based on the reduced slow system, we can conduct  parameter estimation \cite{Ren2015b}, data assimilation \cite{Zhang2017,Qiao2018}, and stochastic bifurcation \cite{Ziying2018}.

\section{Applications}\label{sec3}
In this section, we will discuss some applications associated with the geometric methods for stochastic dynamical systems.

\subsection{Transition phenomena}

In contrast to   deterministic situation, it is   natural to  consider   transitions among  metastable states in stochastic dynamical systems.    That is to say, a transition orbit is likely to  occur  between two stable equilibrium   states, when a system is under random fluctuations.


Transition phenomena occur in gene regulation, tumor cell density, climate change, parametric oscillator, Briggs-Rauscher chemical reaction, and predator-prey systems. By the geometric methods for stochastic dynamical systems, we may qualitatively demonstrate  these noise-induced transition phenomena.

In   genetic regulatory systems, specific protein concentrations play an important role in cell life. The low concentration and high concentration of  specific protein correspond to different cell activities. Recent studies \cite{Xu2013,Zheng2016,Wang2018} have recognized that L\'evy motion can induce switches between different  protein concentrations. Multiple phenotypic states often arise in a single cell with different gene expression states. The transitions between   two phenotypic states due to stochastic  fluctuations have been verified in \cite{Gehao2015}.

In population dynamics, determining the amount of predation which a prey population can sustain without endangering its survival is an important problem \cite{Horsthemke2006}. The problem is to find out the best strategy for a good management of biological resources. We may also think in the opposite direction, driving a prey population to extinction by sufficient predation. In this case the extinction problem does not only find applications in pest control but also in the medical sciences.

\subsection{Stochastic bifurcation}

Although bifurcation studies for deterministic dynamical systems have a long history, the stochastic bifurcation investigation is still in its early stage. One reason for this slow development   in stochastic bifurcation is due to the lack of  appropriate phase portraits, in contrast to deterministic dynamical systems.

Stochastic bifurcations have been observed in a wide range of nonlinear systems in physical science and engineering\cite{Deco,Bashkirtseva,Bogatenko}. Some works about stochastic bifurcation are analytical studies of invariant measures, together with their spectra and supports.   We recently have used most probable phase portraits to detect stochastic bifurcation  \cite{WH}.

\section{Conclusion}\label{sec4}

Stochastic dynamical systems are widely arise as mathematical models  in biology, physics, chemistry and engineering.
In this chapter, we have reviewed several geometric methods for stochastic dynamical systems and their applications. By most probable phase portraits and mean phase portraits, we may better understand stochastic dynamical behaviors such as      transition orbits and stochastic bifurcation.  Based on the reduced systems on slow manifolds, we may estimate parameters,  detect stochastic bifurcation,  and conduct data assimilation on lower dimensional slow systems.


\begin{thebibliography}{99}
%
\bibitem[Arnold (2003)]{Arnold2003}
 Arnold, L. [2003]  \emph{Random Dynamical Systems}, 2nd Ed. (Springer , New York ,USA).

\bibitem[Duan (2015)]{Duan2015}  Duan, J., [2015]  \emph {An Introduction to Stochastic Dynamics}  (Cambridge University Press, New York,USA)

\bibitem[Meiss(2017)]{Meiss2017} Meiss, J.D. [2017],  \emph{Differential Dynamical Systems } Revised. (Society for Industrial and Applied, Philadelphia).

\bibitem[Wiggins (2003)]{Wiggins2003}  Wiggins, S.,  [2003]  \emph{ Introduction to Applied Nonlinear Dynamical Systems and Chaos},  2nd Ed. (Springer , New York ,USA).

\bibitem[Guckenheimer \& Holmes (1983)]{Guckenheimer1983}
 Guckenheimer, J. and   Holmes, P., [1983] \emph{ Nonlinear Oscillations,
Dynamical Systems and Bifurcations of Vector Fields}, (Springer , New York ,USA).

\bibitem[Applebaum (2009)]{Applebaum2009}   Applebaum, D.,  [2009] \emph{ L\'evy Processes and Stochastic Calculus} 2nd Edition. (Cambridge University Press, New York,USA)

\bibitem[Sato (1999)]{Sato1999}  Sato, K.,   [1999]  \emph{ L\'evy Processes and Infinitely Divisible Distributions} (Cambridge University Press, New York,USA)


\bibitem[Holloway \& Spirov (2017)]{Holloway2017}  Holloway, D. M. and  Spirov, A. V., [2017]. Transcriptional bursting inDrosophiladevelopment: Stochastic dynamics ofevestripe 2 expression, \emph{Plos One }  {\bf 12 (4)}, e0176228.

\bibitem[Kumar {\it et al.} (2015)] {Kumar2015} Kumar, N. ,  Singh,  A. and  Kulkarni, R. V., [2015]  Transcriptional bursting in gene expression: analytical results for general stochastic models, \emph{Plos Computational Biology}  {\bf 11}, e1004292.

\bibitem[Dar  {\it et al.} (2012)]{Dar2012} Dar,  R. D.,    Razooky , B. S.,   Singh , A., Trimeloni ,  T. V.,   McCollum , J. M.,
 Cox ,C. D., Simpson,  M. L. and  Weinberger, L. S.,  [2012]  Transcriptional burst frequency and burst size are equally modulated across the human genome, \emph{ PNAS} {\bf 109}, 17454 - 17459.

\bibitem[Ozbudak  {\it et al.} (2002)]{Ozbudak2002}  Ozbudak, E. M.,  Thattai, M.,  Kurtser, I. ,  Grossman, A. D. and  Oudenaarden,  A. V., [2002]. Regulation of noise in the expression of a single gene, \emph{Nature Genetics}  {\bf31}, 69-73.


\bibitem[Blake  {\it et al.} (2003)]{Blake2003}   Blake, W. J.,   K{\ae}rn, M.,  Cantor, C. R.and  Collins, J. J., [2003]. Noise in eukaryotic gene expression, \emph{Nature}  {\bf422}, 633-637.

\bibitem[Sanchez \&  Golding (2013)]{Sanchez2013}  Sanchez, A. and  Golding, I., [2013]. Genetic determinants and cellular constraints in noisy gene expression, \emph{Science} {\bf342(6163)}, 1188-1193.


\bibitem[Karatzas \& Shreve (1991)]{Karatzas1991}  Karatzas,  I. and Shreve, S. E. [1991]  \emph{Brownian Motion and Stochastic Calculus}, 2nd Edition. (Springer-Verlag , New York ,USA)

\bibitem[Arnold (1988)]{Arnold1988}
 Arnold, I. [1988]  \emph{Geometrical Methods in the Theory of Ordinary Differential Equations}, 2nd Ed. (Springer-Verlag , New York ,USA).

\bibitem[Arnold (1989)]{Arnold1989}
 Arnold, L. [1989]  \emph{Mathematical Methods of Classical Mechanics}, 2nd Ed. (Springer-Verlag , New York ,USA).


\bibitem[Cheng { et al.} (2016)]{Cheng2016}   Cheng ,  Z. ,    Duan , J. and  Wang,  L., [2016]  Most probable dynamics of some nonlinear systems under noisy fluctuatuons, \emph{ Commun. Nonlinear Sci. Numer. Simulat.}  {\bf 30}, 108-114.

\bibitem[Duan { et al.} (2003)]{Duan2003} Duan, J.,  Lu, K. and  Schmalfuss, B., [2003]. Invariant manifolds for stochastic partial differential
equations. \emph{ Annals of Probability}.  {\bf 31}, 2109-2135.

\bibitem[ Schmalfuss \&  Schneider(2008)]{Fu2008} Schmalfuss, B. and Schneider, R., [2008]. Invariant manifolds for random dynamical systems with slow and fast
variables. \emph{ Journal of Dynamics and Differential Equations}, {\bf 20}, 133-164.

 \bibitem[Fu { et al.} (2013)]{Fu2013} Fu, H., Liu, X.  and Duan, J., [2013]. Slow manifolds for multi-time-scale stochastic evolutionary systems.  \emph{ Communications
in Mathematical Sciences}. {\bf11(1)}, 141-162.

\bibitem[Ren{ et al.} (2015a)]{Ren2015a}Ren, J., Duan, J. and Christopher, J., [2015]. Approximation of random slow manifolds and settling of inertial particles under uncertainty. \emph{Journal of Dynamics and Differential Equations} {\bf 27},  961-979.

\bibitem[Ren{ et al.} (2015b)]{Ren2015b} Ren J., Duan J. and Wang, X., [2015]. A parameter estimation method based on random slow manifolds. \emph{ Applied Mathematical Modelling} {\bf 39}, 3721-3732.

 \bibitem[Constable { et al.} (2013)]{Constable2013}  Constable, G., McKane, A. and Rogers, T., [2013] Stochastic dynamics on slow manifolds. \emph{Journal of Physics
A: Mathematical and Theoretical}, {\bf 46(295002)}.


\bibitem[Zhang{ et al.}(2017)]{Zhang2017} Zhang, Y., Cheng, Z., Zhang, X., Chen, X., Duan, J. and Li, X. [2017]. Data assimilation and parameter estimation for a multiscale stochastic system with $\alpha$-stable L\'evy noise. \emph{ Journal of Statistical Mechanics: Theory and Experiment}, 113401.


\bibitem[Qiao { et al.} (2018)]{Qiao2018}  Qiao, H., Zhang, Y. and Duan, J. [2018]. Effective filtering on a random slow manifold. \emph{ Nonlinearity} {\bf 31}, 4649¨C4666.

\bibitem[He { et al.} (2018)]{Ziying2018} He, Z., Cai, R., Duan, J. and Liu, X. [2018]. A role of random slow manifolds in detecting stochastic bifurcation. arxiv: 1805.04653v1.

\bibitem[Xu { et al.} (2013)]{Xu2013} Xu, Y., Feng,  J.,  Li, J.and  Zhang,  H.,  [2013] L\'evy noise induced switch in the gene transcriptional regulatory system, \emph{Chaos}  {\bf23(1)}, 013110 .

\bibitem[Zheng { et al.} (2016)]{Zheng2016} Zheng,  Y.,   Larissa,S.,   Duan, J. and   J\"{u}rgen, K.,   [2016]  Transitions in a genetic transcriptional regulatory system under L\'evy motion,  {\it Scient. Rep.} {\bf 6}, 29274.


\bibitem[Wang { et al.} (2018)]{Wang2018}Wang, H., Cheng, X., Duan, J., J\"{u}rgen, K. and Li, X.  [2018] Likelihood for transcriptions in a genetic regulatory system under asymmetric stable L¨¦vy noise. \emph{Chaos},   {\bf28(1)}, 013121.

\bibitem[Ge { et al.} (2015)]{Gehao2015} Ge, H.,  Qian, H. and  Xie, X.S.,  [2015] Stochastic phenotype transition of a single cell in an intermediate region of gene state switching, \emph{Phys. Rev. Lett.}  {\bf114 (7)}, 078101.

\bibitem[Horsthemke \& Lefever(2006)]{Horsthemke2006} Horsthemke, W. and Lefever, R. [2006],   \emph{Noise-Induced Transitions}. (Springer:Berlin Heidelberg ).


\bibitem[Deco \& Mart\'i (2007)]{Deco}  Deco, G.,  Mart\'i, D., [2007],  Deterministic analysis of stochastic bifurcations in multi-stable neurodynamical systems. \emph{Biological Cybernetics},   {\bf96(5)}:487-496.

\bibitem[Bashkirtseva { et al.} (2018)]{Bashkirtseva}    Bashkirtseva, I., Nasyrova, V., Ryashko, L. [2018] Noise-induced bursting and chaos in the two-dimensional Rulkov model. \emph{Chaos Solitons \& Fractals},   {\bf110}: 76-81.

\bibitem[Bogatenko \& Semenov (2018)]{Bogatenko} Bogatenko, T., Semenov, V.  [2018] Coherence resonance in an excitable potential well.  \emph{Physics Letters A }  {\bf382} 2645¨C2649.

\bibitem[Wang { et al.} (2018)]{WH}  Wang, H.,  Chen, X. and Duan, J. [2018] A Stochastic Pitchfork Bifurcation in Most Probable Phase Portraits. \emph{International Journal of Bifurcation and Chaos}, {\bf 28(1)}


\end{thebibliography}
\end{document}